\documentclass{article}
\usepackage[%
journal=JFG,    
lang=british,   
]{ems-journal}

\usepackage{pgfplots}
\usepackage{float}
\usepackage{adjustbox}
\usepackage{esint}
\allowdisplaybreaks[4]
\numberwithin{figure}{section}
\theoremstyle{plain}
\newtheorem{thm}{\protect\theoremname}[section]
\theoremstyle{definition}

\theoremstyle{definition}

\theoremstyle{plain}
\newtheorem{prop}[thm]{\protect\propositionname}
\theoremstyle{plain}
\newtheorem{lem}[thm]{\protect\lemmaname}

\theoremstyle{plain}

\theoremstyle{definition}

\usepackage{tikz}
\usetikzlibrary{shapes,arrows}
\usepackage{verbatim}

\usepackage{amstext}

\newcommand{\R}{\mathbb{R}}
\newcommand{\N}{\mathbb{N}}

\newcommand{\udim}{\overline{\dim}_{\mathrm B}\,}
\newcommand{\ldim}{\underline{\dim}_{\mathrm B}\,}
\newcommand{\bdim}{\dim_{\mathrm B}}

\newcommand{\adim}{\dim_{\mathrm A}}
\newcommand{\ladim}{\dim_{\mathrm L}}

\makeatother

\providecommand{\corollaryname}{Corollary}
\providecommand{\definitionname}{Definition}
\providecommand{\examplename}{Example}
\providecommand{\lemmaname}{Lemma}

\providecommand{\propositionname}{Proposition}
\providecommand{\remarkname}{Remark}
\providecommand{\theoremname}{Theorem}

\numberwithin{equation}{section}

\begin{document}

\title{Assouad type spectra for some self-affine sponges}
\titlemark{Assouad type spectra for some self-affine sponges}



\emsauthor{1}{
	\givenname{Qiang}
	\surname{Huo}
	\mrid{1576706}
	\orcid{0000-0001-9995-4899}}{Q.~Huo}

\Emsaffil{1}{
	\department{School of Mathematical Sciences}
	\organisation{University of Science and Technology of China}
	\rorid{04c4dkn09}
	\address{No. 96, JinZhai Road}
	\zip{230026}
	\city{Hefei}
	\country{China}
	\affemail{qianghuo@ustc.edu.cn}}

\classification[37C45]{28A80}


\begin{abstract}
In this paper we compute the Assouad and lower spectra of Bedford--McMullen sponges in $\R^3$ explicitly. According to the formulae established, we discover that the spectra are not determined by the ratio set, and the box, lower and Assouad dimensions of the sponge anymore, which is unlike the situation in a planar carpet. As a by-product, we construct two Bedford--McMullen sponges on the same grid, both of which have non-uniform fibres. Particularly, they share the same box, lower and Assouad dimensions. However, their Assouad type spectra are different and therefore they are not bi-Lipschitz equivalent. For Bedford--McMullen sponges in higher dimensions, we also determine the dimension spectra when $\theta$ is smaller than the minimal ratio or it is bigger than the maximal ratio.
\end{abstract}

\maketitle


\section{Introduction}
 Let $F\subset\R^d$ be a non-empty compact bounded set. The \textit{Assouad dimension} of $F$ is defined to be the infimum of all exponents $s>0$ for which there exists $C>0$ such that $N(B(x,R)\cap F,r)\leq C(R/r)^s$ for all $x\in F$ and $0<r<R$, where $N(E,r)$ denotes the minimal number of $r$-balls required to cover $E\subset X$ and $B(x,R)$ is the closed ball with centre $x$ and radius $R$. The Assouad dimension, initiated by Assouad \cite{Ass77} in his PhD dissertation, plays a vital role in various fields like embedding theory \cite{Rob09}, number theory \cite{FY18c}, Sobolev inequalities in functional analysis \cite{DLV18,LV16} etc. We refer the readers to books \cite{Fra21,MT10,Rob11} for more background, basic properties and related applications of the Assouad dimension. Unlike the box dimension which captures the global scaling structure of sets, the Assouad dimension quantifies the size of the thickest parts between all pair of scales. To better understand the gap between the box and Assouad dimensions, Fraser and Yu \cite{FY18} defined the \textit{Assouad spectrum} of a set $F$ to be
\begin{align*}
        \adim^{\theta} F=&\inf\left\{s>0: \exists C>0, \text{ s.t. for all } 0<R<1 \text{ and } x\in F, \right. \\
        &\left. N(B(x,R)\cap F,R^{1/\theta})\leq C\left(\frac{R}{R^{1/\theta}}\right)^s\right\}
    \end{align*}
for all $\theta\in(0,1)$. As a dual spectrum, the \textit{lower spectrum} of $F$ is given by
\begin{equation*}
    \begin{split}
        \ladim^{\theta} F=&\sup\left\{s>0: \exists C>0, \text{ s.t. for all } 0<R<1 \text{ and } x\in F, \right. \\
        &\left. N(B(x,R)\cap F,R^{1/\theta})\geq C\left(\frac{R}{R^{1/\theta}}\right)^s\right\}
    \end{split}
    \end{equation*}
for all $\theta\in(0,1)$. In these notions of dimension spectra, the relationship between the pair of scales $R,r=R^{1/\theta}$ is fixed by the parameter $\theta$. The Assouad spectrum interpolates between the well-known upper box dimension and the quasi-Assouad dimension introduced by L\"{u} and Xi \cite{LX16}, i.e. it approaches the upper box dimension as $\theta\to0$ and tends to the quasi-Assouad dimension as $\theta\to1$. Compared with the Assouad and lower dimensions, the corresponding dimension spectra provide us with finer geometric information about the scaling structure of a set. However, the Assouad spectrum is not necessarily monotonic, as shown in \cite[Section 8]{FY18}. 
We mention a variant spectrum introduced in \cite{FHHTY19} and coupled with ``a tale of two spectra'' analysis, called the upper spectrum, which only fixes $R^{1/\theta}$ as an upper bound to $r$ in the original definition of the Assouad spectrum and turns out to be non-decreasing in $\theta$. Beyond the monotonic property, Rutar \cite{Rut24} investigated all possible forms of the Assouad spectrum of a bounded set $F\subset\R^d$.

In \cite{FY18b}, Fraser and Yu obtained the precise formulae for the spectra of planar Bedford--McMullen carpets. They showed that the Assouad spectrum of a Bedford--McMullen carpet with non-uniform fibres has uniquely one phase transition, which occurs when the spectrum approaches the Assouad dimension. Perhaps Bedford--McMullen carpets are the simplest and oldest self-affine sets. Nowadays there are some more general self-affine sets well-studied, such as Lalley--Gatzouras carpets \cite{GL92}, Bara\'{n}ski carpets \cite{Bar07}, and their extension in higher dimensions \cite{How19,KP96} (i.e. the ambient spatial dimension is at least $3$). The precise formulae for Assouad type dimensions of these general self-affine sets can be found in \cite{Fra14,FH17,How19,Mac11}. Recently, Banaji, Fraser, Kolossv\'{a}ry and Rutar \cite{BFKR24} determined the Assouad spectrum of Lalley--Gatzouras carpets and discussed the differentiability of the spectrum. It is natural to study the dimension spectra of \textit{self-affine sponges}, including the Bedford--McMullen class and the Lalley--Gatzouras class.

In this paper, we calculate the precise formulae for the Assouad and lower spectra of Bedford--McMullen sponges in $\R^3$. The key ingredient in our proof is the covering argument of approximate cubes \cite{Ols98}, which follows Fraser and Yu's treatment of planar Bedford--McMullen carpets in \cite{FY18b} essentially. As a consequence, we may construct two sponges in $\R^3$ which are not bi-Lipschitz equivalent. This phenomenon is revealed only by the dimension spectra as these two sponges are defined on the same grid and share the same box, lower and Assouad dimensions. We emphasize that our construction is essentially different from that in \cite[Proposition 3.4]{FY18b} since they considered two carpets defined on grids of size $m\times n$ and $m'\times n'$ respectively and ensured that $\log n/\log m\neq\log n'/\log m'$. 
Fraser and Yu \cite[Corollary 3.5]{FY18b} also proved that the dimension spectra of a planar Bedford--McMullen carpet are completely determined by the ratio $\log n/\log m$ and the box, lower and Assouad dimensions of the carpet. However, this feature is not exhibited by Bedford--McMullen sponges in $\R^3$ anymore according to the explicit formulae established in this paper. Even it is rather complicated to calculate the dimension spectra for sponges in general $\R^d$, inspired by the heuristic covering strategy to deal with approximate cubes in $\R^3$, we may obtain the spectra formulae of Bedford--McMullen sponges in $\R^d$ when $\theta$ is smaller than the minimal ratio or it is bigger than the maximal ratio.

It is worth mentioning that a related notion via dimension interpolation between the Hausdorff and box dimensions, named as the \textit{intermediate dimensions}, was considered in \cite{FFK20}.
Recently, Banaji and Kolossv\'{a}ry \cite{BK24} calculated the intermediate dimensions explicitly for all Bedford--McMullen carpets and revealed the bi-Lipschitz equivalence of two Bedford--McMullen carpets from the perspective of intermediate dimensions. We refer interested readers to \cite{Fra21b} for these two approaches to dimension interpolation.

\section{Preliminaries}

    Throughout the paper, for real-valued functions $f$ and $g$, we write $f\lesssim g$ if there exists $C>0$ such that $f(x)\leq C g(x)$ for all $x$. Sometimes, we use the notation $f(x)=O(g(x))$ to mean the same thing. Similarly, we say that $f\gtrsim g$ if $g\lesssim f$. If both $f\gtrsim g$ and $f\lesssim g$, then we write $f\asymp g$.

    Let $F\subset\R^d$ be a non-empty compact bounded set. The \textit{upper} and \textit{lower box dimensions} of $F$ are defined by
    \begin{equation*}
        \udim F=\limsup_{r\to0}\frac{\log N(F,r)}{\log (1/r)}  \text{  and  } \ldim F=\liminf_{r\to0}\frac{\log N(F,r)}{\log (1/r)}.
    \end{equation*}
    If these two values coincide, the common value is called the \textit{box dimension} of $F$ and denoted by $\bdim F$. Following \cite{Fra14}, we define the \textit{lower dimension} of $F$ by
    \begin{align*}
        \ladim F=&\sup\left\{s>0: \exists C>0, \text{ s.t. for all } 0<r<R \text{ and } x\in F, \right. \\
        &\left. N(B(x,R)\cap F,r)\geq C\left(\frac{R}{r}\right)^s\right\}.
    \end{align*}
   In particular, Fraser and Yu \cite[Propositions 3.1 and 3.9]{FY18} proved the general bound on the Assouad type spectra in terms of the box, lower and Assouad dimensions as follows. For any $\theta\in(0,1)$ one has
\begin{equation*}
    \udim F\leq\adim^{\theta} F\leq\min\left\{\frac{\udim F}{1-\theta}, \adim F\right\} \text{ and } \ladim F\leq \ladim^{\theta} F\leq \ldim F.
\end{equation*}

\section{Bedford--McMullen sponges and main results}
We begin with the notion of Bedford--McMullen sponges initiated by Kenyon and Peres \cite{KP96}.
Fix $d\in\N$ and let $2\leq n_1<\ldots<n_d$ be integers. Let $\mathcal{I}_q=\{0,\ldots,n_q-1\}$ and choose a subset $\mathcal{D}\subset \prod_{q=1}^d\mathcal{I}_q$ composed of at least two elements. For each $\bar\imath=(i_1,\ldots,i_d)\in\mathcal{D}$, we define an affine transform $\varphi_{\bar\imath}:[0,1]^d\to[0,1]^d$ by
\begin{equation*}
    \varphi_{\bar\imath}(x_1,\ldots,x_d)=\left(\frac{x_1+i_1}{n_1},\ldots,\frac{x_d+i_d}{n_d}\right).
\end{equation*}
Let $\mathcal{D}^{\infty}=\{\omega=(\bar\imath_1,\bar\imath_2,\ldots):\bar\imath_l=(i_{l,1},\ldots,i_{l,d})\in\mathcal{D}\}$ be the set of all infinite words over $\mathcal{D}$ and let $\Pi:\mathcal{D}^{\infty}\to[0,1]^d$ be the canonical map from the symbolic space to the geometric space defined by
\begin{equation*}
    \Pi(\omega)=\bigcap\limits_{l\in\N}\varphi_{\bar\imath_1}\circ\cdots\circ\varphi_{\bar\imath_l}([0,1]^d).
\end{equation*}
The set $F:=\Pi(\mathcal{D})$ is called the \textit{Bedford--McMullen sponge}. Denote the \textit{ratio set} of $F$ by $\{\log n_p/\log n_q\}_{1\leq p<q\leq d}$.
For $l=1,2,\ldots,d$, we define $\pi_{l}:\mathcal{D}\to\prod_{k=1}^{l}\mathcal{I}_k$ to be the projection onto the first $l$ coordinates. For $l=1,2,\ldots,d-1$ and $(i_1,\ldots,i_l)\in\pi_{l}\mathcal{D}$, let
\begin{equation*}
     N(i_1,\ldots,i_l)=|\{i_{l+1}\in\mathcal{I}_{l+1}:(i_1,\ldots,i_l,i_{l+1})\in\pi_{l+1}\mathcal{D}\}|.
\end{equation*}
We say that a Bedford--McMullen sponge has \textit{uniform fibres} if and only if for all $l=1,2,\ldots,d-1$ it holds
\begin{equation*}
    N(i_1,\ldots,i_l)=N(j_1,\ldots,j_l)
\end{equation*}
for all $(i_1,\ldots,i_l),(j_1,\ldots,j_l)\in\pi_{l}\mathcal{D}$.

Formulae for the Hausdorff, box, Assouad and lower dimensions of Bedford--McMullen sponges were established in \cite{FH17,KP96}, which can be seen as generalisation of the formulae in the planar case obtained by Bedford \cite{Bed84}, McMullen \cite{Mc84}, Mackay \cite{Mac11} and Fraser \cite{Fra14}.
\begin{thm}[Kenyon--Peres--Fraser--Howroyd{}]\label{Kenyon-Peres-Fraser-Howroyd}
    The Assouad, box and lower dimensions of a Bedford--McMullen sponge $F$ are given by
    \begin{equation*}
    \begin{split}
        &\adim F=\frac{\log |\pi_1\mathcal{D}|}{\log n_1}+\sum\limits_{l=1}^{d-1}\frac{\log \max\limits_{(i_1,\ldots,i_l)\in\pi_{l}\mathcal{D}}N(i_1,\ldots,i_l)}{\log n_{l+1}},\\
        &\bdim F=\frac{\log|\pi_1\mathcal{D}|}{\log n_1}+\sum\limits_{l=1}^{d-1}\frac{\log (|\pi_{l+1}\mathcal{D}|/|\pi_{l}\mathcal{D}|)}{\log n_{l+1}},\\
        &\ladim F=\frac{\log |\pi_1\mathcal{D}|}{\log n_1}+\sum\limits_{l=1}^{d-1}\frac{\log \min\limits_{(i_1,\ldots,i_l)\in\pi_{l}\mathcal{D}}N(i_1,\ldots,i_l)}{\log n_{l+1}}.
    \end{split}
    \end{equation*}
    Moreover, either $\ladim F<\bdim F<\adim F$ (the `non-uniform fibres' case) or $\ladim F=\bdim F=\adim F$ (the `uniform fibres' case).
\end{thm}

\subsection*{Notation 1}
Let $F$ be a Bedford--McMullen sponge in $\R^3$ with non-uniform fibres. Denote by 
\begin{equation*}
\begin{split}
    & M_{\max}=\max_{i_1\in\pi_1 \mathcal{D}} N(i_1), M_{\min}=\min_{i_1\in\pi_1 \mathcal{D}} N(i_1);\\
    & N_{\max}=\max_{(i_1,i_2)\in\pi_2 \mathcal{D}}N(i_1,i_2), N_{\min}=\min_{(i_1,i_2)\in\pi_2 \mathcal{D}}N(i_1,i_2);\\
    & C_{\max}=\max_{i_1\in\pi_1 \mathcal{D}}\sum_{i_2: (i_1,i_2)\in\pi_2 \mathcal{D}} N(i_1,i_2), C_{\min}=\min_{i_1\in\pi_1 \mathcal{D}}\sum_{i_2: (i_1,i_2)\in\pi_2 \mathcal{D}} N(i_1,i_2).
\end{split}
\end{equation*} 
Using these notations, Theorem \ref{Kenyon-Peres-Fraser-Howroyd} yields that
 \begin{equation*}
 \begin{split}
     &\adim F=\frac{\log |\pi_1\mathcal{D}|}{\log n_1}+\frac{\log M_{\max}}{\log n_2}+\frac{\log N_{\max}}{\log n_3},\\
     &\bdim F=\frac{\log|\pi_1\mathcal{D}|}{\log n_1}+\frac{\log (|\pi_2\mathcal{D}|/|\pi_1\mathcal{D}|)}{\log n_2}+\frac{\log (|\mathcal{D}|/|\pi_2\mathcal{D}|)}{\log n_3},\\
     &\ladim F=\frac{\log |\pi_1\mathcal{D}|}{\log n_1}+\frac{\log M_{\min}}{\log n_2}+\frac{\log N_{\min}}{\log n_3}.
 \end{split}
 \end{equation*}

For $\theta\in(0,1)$, denote by
\begin{align*}
        f_1(\theta)&=\frac{\bdim F-\theta\left(\frac{\log (|\mathcal{D}|/C_{\max})}{\log n_1}+\frac{\log (C_{\max}/N_{\max})}{\log n_2}+\frac{\log N_{\max}}{\log n_3}\right)}{1-\theta},\\
        f_2(\theta)&=\frac{\frac{\log|\pi_1\mathcal{D}|}{\log n_1}+\frac{\log (|\pi_2\mathcal{D}|/|\pi_1\mathcal{D}|)}{\log n_2}+\frac{\log (C_{\max}/M_{\max})}{\log n_3}}{1-\theta}\\
        &-\frac{\theta\left(\frac{\log (|\pi_2\mathcal{D}|/M_{\max})}{\log n_1}+\frac{\log (C_{\max}/N_{\max})}{\log n_2}+\frac{\log N_{\max}}{\log n_3}\right)}{1-\theta},\\
        f_3(\theta)&=\frac{\log |\pi_1\mathcal{D}|}{\log n_1}+\frac{\frac{\log M_{\max}}{\log n_2}+\frac{\log (C_{\max}/M_{\max})}{\log n_3}-\theta\left(\frac{\log (C_{\max}/N_{\max})}{\log n_2}+\frac{\log N_{\max}}{\log n_3}\right)}{1-\theta},\\
        f_4(\theta)&=\frac{\frac{\log|\pi_1\mathcal{D}|}{\log n_1}+\frac{\log (|\pi_2\mathcal{D}|/|\pi_1\mathcal{D}|)}{\log n_2}-\theta\left(\frac{\log (|\pi_2\mathcal{D}|/M_{\max})}{\log n_1}+\frac{\log M_{\max}}{\log n_2}\right)}{1-\theta}+\frac{\log N_{\max}}{\log n_3}
\end{align*}
and
\begin{align*}
        g_1(\theta)&=\frac{\bdim F-\theta\left(\frac{\log (|\mathcal{D}|/C_{\min})}{\log n_1}+\frac{\log (C_{\min}/N_{\min})}{\log n_2}+\frac{\log N_{\min}}{\log n_3}\right)}{1-\theta},\\
        g_2(\theta)&=\frac{\frac{\log|\pi_1\mathcal{D}|}{\log n_1}+\frac{\log (|\pi_2\mathcal{D}|/|\pi_1\mathcal{D}|)}{\log n_2}+\frac{\log (C_{\min}/M_{\min})}{\log n_3}}{1-\theta}\\
        &-\frac{\theta\left(\frac{\log (|\pi_2\mathcal{D}|/M_{\min})}{\log n_1}+\frac{\log (C_{\min}/N_{\min})}{\log n_2}+\frac{\log N_{\min}}{\log n_3}\right)}{1-\theta},\\
        g_3(\theta)&=\frac{\log |\pi_1\mathcal{D}|}{\log n_1}+\frac{\frac{\log M_{\min}}{\log n_2}+\frac{\log (C_{\min}/M_{\min})}{\log n_3}-\theta\left(\frac{\log (C_{\min}/N_{\min})}{\log n_2}+\frac{\log N_{\min}}{\log n_3}\right)}{1-\theta},\\
        g_4(\theta)&=\frac{\frac{\log|\pi_1\mathcal{D}|}{\log n_1}+\frac{\log (|\pi_2\mathcal{D}|/|\pi_1\mathcal{D}|)}{\log n_2}-\theta\left(\frac{\log (|\pi_2\mathcal{D}|/M_{\min})}{\log n_1}+\frac{\log M_{\min}}{\log n_2}\right)}{1-\theta}+\frac{\log N_{\min}}{\log n_3}.
\end{align*}

We are now in position to state the main result (including three pieces), which gives the precise formulae for the Assouad and lower spectra of a Bedford--McMullen sponge in $\R^3$. We always assume that the Bedford--McMullen sponge has non-uniform fibres.

\begin{thm}\label{spectra formula}
    Let $F$ be a Bedford--McMullen sponge in $\R^3$ satisfying $\log n_1/\log n_2<\log n_2/\log n_3$. Then the Assouad spectrum of $F$ is given by
    \begin{equation*}
      \adim^{\theta}F=  
       \begin{cases}
         f_1(\theta) & \text{for }\, 0<\theta\leq\frac{\log n_1}{\log n_3},\\
         f_2(\theta) & \text{for }\, \frac{\log n_1}{\log n_3}<\theta\leq\frac{\log n_1}{\log n_2},\\
         f_3(\theta) & \text{for }\, \frac{\log n_1}{\log n_2}<\theta\leq\frac{\log n_2}{\log n_3},\\
         \adim F   &   \text{for }\, \frac{\log n_2}{\log n_3}<\theta<1
       \end{cases}
    \end{equation*}
    and the lower spectrum of $F$ is given by
    \begin{equation*}
      \ladim^{\theta}F=  
       \begin{cases}
         g_1(\theta) & \text{for }\, 0<\theta\leq\frac{\log n_1}{\log n_3},\\
         g_2(\theta) & \text{for }\, \frac{\log n_1}{\log n_3}<\theta\leq\frac{\log n_1}{\log n_2},\\
         g_3(\theta) & \text{for }\, \frac{\log n_1}{\log n_2}<\theta\leq\frac{\log n_2}{\log n_3},\\
         \ladim F   &   \text{for }\, \frac{\log n_2}{\log n_3}<\theta<1.
       \end{cases}
    \end{equation*}

\end{thm}

\begin{thm}\label{spectra case(2)}
    Let $F$ be a Bedford--McMullen sponge in $\R^3$ satisfying $\log n_2/\log n_3<\log n_1/\log n_2$. Then the Assouad spectrum of $F$ is given by
        \begin{align*}
      \adim^{\theta}F=  
       \begin{cases}
         f_1(\theta) & \text{for }\, 0<\theta\leq\frac{\log n_1}{\log n_3},\\
         f_2(\theta) & \text{for }\, \frac{\log n_1}{\log n_3}<\theta\leq\frac{\log n_2}{\log n_3},\\
         f_4(\theta) & \text{for }\, \frac{\log n_2}{\log n_3}<\theta\leq\frac{\log n_1}{\log n_2},\\
         \adim F   &   \text{for }\, \frac{\log n_1}{\log n_2}<\theta<1
       \end{cases}
    \end{align*}
    and the lower spectrum of $F$ is given by
    \begin{align*}
      \ladim^{\theta}F=  
       \begin{cases}
         g_1(\theta) & \text{for }\, 0<\theta\leq\frac{\log n_1}{\log n_3},\\
         g_2(\theta) & \text{for }\, \frac{\log n_1}{\log n_3}<\theta\leq\frac{\log n_2}{\log n_3},\\
         g_4(\theta) & \text{for }\, \frac{\log n_2}{\log n_3}<\theta\leq\frac{\log n_1}{\log n_2},\\
         \ladim F   &   \text{for }\, \frac{\log n_1}{\log n_2}<\theta<1.
       \end{cases}
    \end{align*}

\end{thm}

\begin{thm}\label{spectra case(3)}
    Let $F$ be a Bedford--McMullen sponge in $\R^3$ satisfying $\log n_1/\log n_2=\log n_2/\log n_3$. Then the Assouad spectrum of $F$ is given by
    \begin{equation*}
      \adim^{\theta}F=  
       \begin{cases}
         f_1(\theta) & \text{for }\, 0<\theta\leq\frac{\log n_1}{\log n_3},\\
         f_2(\theta) & \text{for }\, \frac{\log n_1}{\log n_3}<\theta\leq\frac{\log n_2}{\log n_3},\\
         \adim F   &   \text{for }\, \frac{\log n_2}{\log n_3}<\theta<1
       \end{cases}
    \end{equation*}
    and the lower spectrum of $F$ is given by
    \begin{equation*}
      \ladim^{\theta}F=  
       \begin{cases}
         g_1(\theta) & \text{for }\, 0<\theta\leq\frac{\log n_1}{\log n_3},\\
         g_2(\theta) & \text{for }\, \frac{\log n_1}{\log n_3}<\theta\leq\frac{\log n_2}{\log n_3},\\
         \ladim F   &   \text{for }\, \frac{\log n_2}{\log n_3}<\theta<1.
       \end{cases}
    \end{equation*}

\end{thm}

 Let $F$ be a Bedford--McMullen carpet with non-uniform fibres defined on a $m\times n$ grid. The ratio $\log m/\log n$ turns out to be a bi-Lipschitz invariant, which was proven by Fraser and Yu in \cite[Theorem 3.3]{FY18b}. Moreover, they \cite[Corollary 3.5]{FY18b} observed that the Assouad and lower spectra of such a carpet are completely determined by the ratio $\log m/\log n$ and the box, lower and Assouad dimensions of the carpet. However, the situation will be very different for Bedford--McMullen sponges. Observe from Theorem \ref{Kenyon-Peres-Fraser-Howroyd} that, for the Assouad dimension, the term $\log \max\limits_{(i_1,\ldots,i_l)\in\pi_{l}\mathcal{D}}N(i_1,\ldots,i_l)/\log n_{l+1}$ relates to the maximal fibre dimension regarding the first $l$ coordinates of $\bar\imath\in\mathcal{D}$; and for the box dimension, we may understand the term $\log (|\pi_{l+1}\mathcal{D}|/|\pi_{l}\mathcal{D}|)/\log n_{l+1}$ as the average fibre dimension in coordinate $l$ by regarding $|\pi_{l+1}\mathcal{D}|/|\pi_{l}\mathcal{D}|$ as the arithmetic average of the $N(i_1,\ldots,i_l)$. 
 
 When considering (for instance) the Assouad spectrum of a Bedford--McMullen sponge defined on a $n_1\times n_2\times n_3$ grid, three new terms $\log(|\mathcal{D}|/C_{\max}), \log(C_{\max}/N_{\max})$, $\log(C_{\max}/M_{\max})$ involving $C_{\max}$ appear, as shown in Theorems \ref{spectra formula}-\ref{spectra case(3)}. To understand $C_{\max}$, fix $i_1\in\pi_1\mathcal{D}$, we group together level-$2$ fibres over the first coordinate $i_1$ and count the total number of next level-$3$ fibres. Then we choose $i_1\in\pi_1\mathcal{D}$ which maximizes the quantity $\sum_{i_2: (i_1,i_2)\in\pi_2\mathcal{D}}N(i_1,i_2)$ (depending only on the level-$1$ fibre). Since these terms relating to $C_{\max}$ or $C_{\min}$ do not appear in the formulae of $\adim F,\bdim F$ and $\ladim F$, the dimension spectra of a sponge are not only determined by the ratio set $\{\log n_p/\log n_q\}_{p<q}$ and the box, lower, Assouad dimension of that sponge anymore. This is the main point different from the planar case. Furthermore, having the same ratio set $\{\log n_p/\log n_q\}_{p<q}$ is a necessary condition for two Bedford--McMullen sponges with non-uniform fibres to be bi-Lipschitz equivalent, but it is not sufficient. Indeed, for two sponges defined on the same grid, we may guarantee the spectra would be different by ensuring that $C_{\max}\neq C'_{\max}$ and $C_{\min}\neq C'_{\min}$ whereas the other parameters in the formulae of the spectra coincide, which are necessary in the construction.
In the following example, we construct two Bedford--McMullen sponges which are not bi-Lipschitz equivalent. This is revealed by the dimension spectra, not by any of the aforementioned dimensions nor the ratio set.
\begin{prop}\label{Lipschitz equivalence}
    We can find two topologically equivalent Bedford--McMullen sponges on the same grid, for which the box, lower and Assouad dimensions coincide. However, their Assouad and lower spectra are different and therefore they are not bi-Lipschitz equivalent.
\end{prop}

After we give the formulae for the dimension spectra of Bedford--McMullen sponges in $\R^3$, it is natural to consider higher-dimensional sponges. For a sponge defined on a $n_1\times n_2\times\cdots\times n_d$ grid, the ratio set is composed of at most $d(d-1)/2$ points in terms of $\log n_p/\log n_q$ with $1\leq p<q\leq d$. Following the proof of Theorem \ref{spectra formula}, we use different covering strategy to deal with approximate cubes when $\theta<\log n_p/\log n_q$ and $\theta>\log n_p/\log n_q$, which will possibly lead to different formulae for the spectra. This means that the dimension spectra probably fail to be differentiable at these points which belong to the ratio set.
Moreover, we need to consider the generalization of $C_{\max}$ and $C_{\min}$ to higher dimensions, which will make the calculation rather complicated. For these reasons, we do not expect to compute the spectra $\adim^{\theta}F,\ladim^{\theta}F$ for all $\theta\in(0,1)$. Instead, owning to the heuristic proof of Theorem \ref{spectra formula}, we discover that there are two situations for which we can easily obtain the formulae for the spectra, i.e. when $\theta$ is smaller than the minimal ratio or $\theta$ is bigger than the maximal ratio and smaller than $1$.

\subsection*{Notation 2}
Let $F$ be a Bedford--McMullen sponge in $\R^d$ with non-uniform fibres. Denote by 
\begin{align*}
     N_{d-1,\max}&=\max_{(i_1,\ldots,i_{d-1})\in\pi_{d-1} \mathcal{D}}N(i_1,\ldots,i_{d-1}),\\
     N_{d-1,\min}&=\min_{(i_1,\ldots,i_{d-1})\in\pi_{d-1} \mathcal{D}}N(i_1,\ldots,i_{d-1});\\
     N_{k,\max}&=\max_{(i_1,\ldots,i_{k})\in\pi_{k} \mathcal{D}} \sum_{i_{k+1}:(i_1,\ldots,i_{k+1})\in\pi_{k+1}\mathcal{D}} \sum_{i_{k+2}:(i_1,\ldots,i_{k+2})\in\pi_{k+2}\mathcal{D}}\\
     &\cdots\sum_{i_{d-1}:(i_1,\ldots,i_{d-1})\in\pi_{d-1}\mathcal{D}} N(i_1,\ldots,i_{d-1}), \\
     N_{k,\min}&=\min_{(i_1,\ldots,i_{k})\in\pi_{k} \mathcal{D}} \sum_{i_{k+1}:(i_1,\ldots,i_{k+1})\in\pi_{k+1}\mathcal{D}} \sum_{i_{k+2}:(i_1,\ldots,i_{k+2})\in\pi_{k+2}\mathcal{D}}\\
     &\cdots\sum_{i_{d-1}:(i_1,\ldots,i_{d-1})\in\pi_{d-1}\mathcal{D}} N(i_1,\ldots,i_{d-1});\\
     M_{k,\max}&=\max_{(i_1,\ldots,i_k)\in\pi_k \mathcal{D}}N(i_1,\ldots,i_k),\\
     M_{k,\min}&=\min_{(i_1,\ldots,i_k)\in\pi_k \mathcal{D}}N(i_1,\ldots,i_k)
\end{align*}
for $k=1,\ldots,d-2$. Then Theorem \ref{Kenyon-Peres-Fraser-Howroyd} yields that
\begin{equation*}
\begin{split}
    &\adim F=\frac{\log |\pi_1\mathcal{D}|}{\log n_1}+\sum\limits_{k=1}^{d-2}\frac{\log M_{k,\max}}{\log n_{k+1}}+\frac{\log N_{d-1,\max}}{\log n_{d}},\\
    &\ladim F=\frac{\log |\pi_1\mathcal{D}|}{\log n_1}+\sum\limits_{k=1}^{d-2}\frac{\log M_{k,\min}}{\log n_{k+1}}+\frac{\log N_{d-1,\min}}{\log n_{d}}.
\end{split}
\end{equation*}

\begin{thm}\label{general dimension}
    Let $F$ be a Bedford--McMullen sponge in $\R^d$. 
    \begin{enumerate}[$(1)$]
        \item If $\theta\in(0,\log n_1/\log n_d]$, then
        \begin{equation*}
           \adim^{\theta}F=\frac{\bdim F-\theta\left(\frac{\log (|\mathcal{D}|/N_{1,\max})}{\log n_1}+\sum\limits_{k=1}^{d-2}\frac{\log (N_{k,\max}/N_{k+1,\max})}{\log n_{k+1}}+\frac{\log N_{d-1,\max}}{\log n_d}\right)}{1-\theta}
       \end{equation*}
       and
       \begin{equation*}
           \ladim^{\theta}F=\frac{\bdim F-\theta\left(\frac{\log (|\mathcal{D}|/N_{1,\min})}{\log n_1}+\sum\limits_{k=1}^{d-2}\frac{\log (N_{k,\min}/N_{k+1,\min})}{\log n_{k+1}}+\frac{\log N_{d-1,\min}}{\log n_d}\right)}{1-\theta};
       \end{equation*}
       \item If $\theta\in[\max_{1\leq p<q\leq d}\log n_p/\log n_q,1)$, then $\adim^{\theta}F=\adim F$ and $\ladim^{\theta}F=\ladim F$.
    \end{enumerate}
\end{thm}

\section{Proofs}
In this section, we prove our main results, Theorems \ref{spectra formula}-\ref{spectra case(3)}, \ref{general dimension} and Proposition \ref{Lipschitz equivalence}.

Fix $\omega=(\bar\imath_1,\ldots)=((i_{1,1},i_{1,2},i_{1,3}),\ldots)\in\mathcal{D}^{\infty}$ and $r>0$ sufficiently small. We denote the \textit{approximate cube} centred at $\omega$ with side length $r>0$ by $Q(\omega,r)$, and define it by
\begin{equation*}
    Q(\omega,r)=\{\omega'=(\bar\jmath_1,\ldots)\in\mathcal{D}^{\infty}: j_{t,q}=i_{t,q}, \forall q=1,2,3 \text{ and } \forall t=1,\ldots,l_q(r)\},
\end{equation*}
where $l_q(r)$ is the unique integer satisfying
\begin{equation}\label{lq}
    n_q^{-l_q(r)} \leq r < n_q^{-l_q(r)+1}.
\end{equation}
In particular,
\begin{equation*}
    -\frac{\log r}{\log n_q} \leq l_q(r) < -\frac{\log r}{\log n_q}+1.
\end{equation*}
The set $Q(\omega,r)$ depends only on the coordinates $i_{1,1},\ldots,i_{l_1(r),1},i_{1,2},\ldots,i_{l_2(r),2},i_{1,3},\ldots$, $i_{l_3(r),3}$. Sometimes we also denote it by
\begin{equation*}
    Q(i_{1,1},\ldots,i_{l_1(r),1},i_{1,2},\ldots,i_{l_2(r),2},i_{1,3},\ldots,i_{l_3(r),3}).
\end{equation*}
The geometric projection of this set, $\Pi(Q(\omega,r))$, is a subset of $F$ which contains $\Pi(\omega)$ and sits inside a cuboid in $\R^3$ which is approximately a cube in the sense that it has a base with length $n_1^{-l_1(r)}\in(r/n_1,r]$ and width $n_2^{-l_2(r)}\in(r/n_2,r]$, and height $n_3^{-l_3(r)}\in(r/n_3,r]$. In order to calculate the Assouad and lower spectra, one may replace $B(x,R)$ by $\Pi(Q(\omega,R))$ since these two quantities give us equivalent definitions for the spectra.

The following lemma is a prior indication that there will be phase transitions in the spectra at $\theta=\log n_1/\log n_3,\log n_1/\log n_2$ and $\log n_2/\log n_3$ (possibly distinct).

\begin{lem}\label{scale comparison}
    Let $R\in(0,1)$ and $\theta\in(0,1)$. We list all possible relationships between $l_q(R)$ and $l_q(R^{1/\theta}), q=1,2,3$ as follows:
    \begin{enumerate}[$(a)$]
        \item $l_3(R)\leq l_2(R) \leq l_1(R) \leq l_3(R^{1/\theta})\leq l_2(R^{1/\theta})\leq l_1(R^{1/\theta})$;
        \item $l_3(R) \leq l_2(R) \leq l_3(R^{1/\theta}) \leq l_1(R) \leq l_2(R^{1/\theta}) \leq l_1(R^{1/\theta})$;
        \item $l_3(R) \leq l_2(R) \leq l_3(R^{1/\theta}) \leq l_2(R^{1/\theta})  \leq l_1(R) \leq l_1(R^{1/\theta})$;
        \item $l_3(R) \leq l_3(R^{1/\theta}) \leq l_2(R) \leq l_2(R^{1/\theta}) \leq l_1(R) \leq l_1(R^{1/\theta})$;
        \item $l_3(R)\leq l_3(R^{1/\theta}) \leq l_2(R)\leq l_1(R)  \leq l_2(R^{1/\theta}) \leq l_1(R^{1/\theta})$.
    \end{enumerate}
    
    \textit{Case 1.} $\log n_1/\log n_2<\log n_2/\log n_3$:
    \begin{enumerate}[$(1)$]
       \item if $\theta\in(0,\log n_1/\log n_3]$, then $(a)$ holds;
       \item if $\theta\in(\log n_1/\log n_3,\log n_1/\log n_2]$, then $(b)$ holds;
       \item if $\theta\in(\log n_1/\log n_2,\log n_2/\log n_3]$, then $(c)$ holds;
       \item if $\theta\in(\log n_2/\log n_3,1)$, then $(d)$ holds.
    \end{enumerate}

    \textit{Case 2.} $\log n_2/\log n_3<\log n_1/\log n_2$:
    \begin{enumerate}[$(1')$]
       \item if $\theta\in(0,\log n_1/\log n_3]$, then $(a)$ holds;
       \item if $\theta\in(\log n_1/\log n_3,\log n_2/\log n_3]$, then $(b)$ holds;
       \item if $\theta\in(\log n_2/\log n_3,\log n_1/\log n_2]$, then $(e)$ holds;
       \item if $\theta\in(\log n_1/\log n_2,1)$, then $(d)$ holds.
    \end{enumerate}

    \textit{Case 3.} $\log n_1/\log n_2=\log n_2/\log n_3$:
    \begin{enumerate}[$(1")$]
       \item if $\theta\in(0,\log n_1/\log n_3]$, then $(a)$ holds;
       \item if $\theta\in(\log n_1/\log n_3,\log n_2/\log n_3]$, then $(b)$ holds;
       \item if $\theta\in(\log n_2/\log n_3,1)$, then $(d)$ holds.
    \end{enumerate}
\end{lem}
\begin{proof}
    It follows directly from the definitions of $l_q(\cdot),q=1,2,3$.
\end{proof}

\subsection{Proof of Theorem \ref{spectra formula}}
    Recall that we are assuming $\log n_1/\log n_2<\log n_2/\log n_3$. Fix $\omega=(\bar\imath_1,\bar\imath_2,\ldots)=((i_{1,1},i_{1,2}$, $i_{1,3}),(i_{2,1},i_{2,2},i_{2,3}),\ldots)\in\mathcal{D}^{\infty}$ and $R\in(0,1)$.
    \begin{enumerate}[$(1)$]
       \item if $\theta\in(0,\log n_1/\log n_3]$:
       
       \textbf{Step $1^{\circ}$:} In order to estimate $N(\Pi(Q(\omega,R)),R^{1/\theta})$, we decompose the approximate cube $Q(i_{1,1},\ldots,i_{l_1(R),1},i_{1,2},\ldots,i_{l_2(R),2},i_{1,3},\ldots,i_{l_3(R),3}) (=Q(\omega,R))$ into cylinders in terms of
       \begin{equation*}
           \{\omega'=(\bar\jmath_1,\ldots)\in Q(\omega,R): j_{l,3}=u_{l,3}, \forall l=l_3(R)+1,\ldots,l_2(R)\},
       \end{equation*}
       where $u_{l,3}\in\mathcal{I}_3$ satisfies that $(i_{l,1},i_{l,2},u_{l,3})\in\mathcal{D}$ for each $l=l_3(R)+1,\ldots,l_2(R)$. The image of each cylinder under the canonical map $\Pi$ sits inside a middle sized cuboid with base length $n_1^{-l_1(R)}$, base width $n_2^{-l_2(R)}$ (which are the same as those of the approximate cube), and height $n_3^{-l_2(R)}$. Moreover, the total number of these cylinders is exactly
       \begin{equation*}
           \prod\limits_{l=l_3(R)+1}^{l_2(R)} N(i_{l,1},i_{l,2}).
       \end{equation*}

       \textbf{Step $2^{\circ}$:} We continue to cover each cylinder in \textbf{Step $1^{\circ}$} independently.
       Fix $(u_{l_3(R)+1,3}, \ldots, u_{l_2(R),3})$ $\in\mathcal{I}_3^{l_2(R)-l_3(R)}$ with $(i_{l,1},i_{l,2},u_{l,3})\in\mathcal{D}$ for each $l=l_3(R)+1,\ldots,l_2(R)$. Denote by 
       \begin{equation*}
           \begin{split}
               &Q(i_{1,1},\ldots,i_{l_1(R),1},i_{1,2},\ldots,i_{l_2(R),2},i_{1,3},\ldots,i_{l_3(R),3}, u_{l_3(R)+1,3},\ldots,u_{l_2(R),3})\\
               &:=\{\omega'=(\bar\jmath_1,\ldots)\in Q(\omega,R): j_{l,3}=u_{l,3}, \forall l=l_3(R)+1,\ldots,l_2(R)\}
           \end{split}
       \end{equation*}
       the cylinder in \textbf{Step $1^{\circ}$}. Observe that this set is further composed of smaller cylinders in terms of \begin{equation*}
         \left\{ \omega'=(\bar\jmath_1,\ldots)\in Q(\omega,R)  \middle|\, 
         \parbox{2.5in}{\centering $j_{l,3}=u_{l,3}, \forall l=l_3(R)+1,\ldots,l_2(R)$ \\
         $j_{l,2}=v_{l,2}$ and $j_{l,3}=v_{l,3}$,\\ $\forall l=l_2(R)+1,\ldots,l_1(R)$ } \right\}, 
     \end{equation*}
     where $(v_{l,2},v_{l,3})\in\mathcal{I}_2 \times \mathcal{I}_3$ satisfies that $(i_{l,1},v_{l,2},v_{l,3})\in\mathcal{D}$ for each $l=l_2(R)+1,\ldots,l_1(R)$. The image of each cylinder under the canonical map $\Pi$ lies inside a smaller sized cuboid with base length $n_1^{-l_1(R)}$, base width $n_2^{-l_1(R)}$ and height $n_3^{-l_1(R)}$. Furthermore, we find exactly  
     \begin{equation*}
           \prod\limits_{l=l_2(R)+1}^{l_1(R)} \sum_{v_{l,2}: (i_{l,1},v_{l,2})\in\pi_2 \mathcal{D}} N(i_{l,1},v_{l,2}).
       \end{equation*}
       cylinders inside $Q(i_{1,1},\ldots,i_{l_1(R),1},i_{1,2},\ldots,i_{l_2(R),2},i_{1,3},\ldots,i_{l_3(R),3}, u_{l_3(R)+1,3}$, $\ldots,u_{l_2(R),3})$. Indeed, for each $l$, the sum over $v_{l,2}$ means that we consider all possible choices of level-$2$ fibres and count the total number of the next level-$3$ fibres thereafter. The quantity $\sum_{v_{l,2}: (i_{l,1},v_{l,2})\in\pi_2 \mathcal{D}} N(i_{l,1},v_{l,2})$ depends only on the level-$1$ fibre $i_{l,1}$ for each $l$.

       \textbf{Step $3^{\circ}$:} After \textbf{Step $2^{\circ}$}, the approximate cube $\Pi(Q(\omega,R))$ is covered by 
       \begin{equation*}
           \prod\limits_{l=l_3(R)+1}^{l_2(R)} N(i_{l,1},i_{l,2})\cdot\prod\limits_{l=l_2(R)+1}^{l_1(R)} \sum_{v_{l,2}: (i_{l,1},v_{l,2})\in\pi_2 \mathcal{D}} N(i_{l,1},v_{l,2})
       \end{equation*}
       cuboids. The height of such cuboid, $n_3^{-l_1(R)}$, is still larger than $R^{1/\theta}/n_3$ since $l_1(R) \leq l_3(R^{1/\theta})$. We want to cover these cuboids by sets of diameter $R^{1/\theta}$ and therefore we iterate the construction inside each cuboid until the height is $n_3^{-l_3(R^{1/\theta})}$. This takes $l_3(R^{1/\theta})-l_1(R)$ iterations and we pick up $|\mathcal{D}|$ smaller cuboids for each iteration.

       \textbf{Step $4^{\circ}$:} After \textbf{Step $3^{\circ}$}, we are left with a large family of cuboids each with height approximately $R^{1/\theta}$, whereas its base length $n_1^{-l_3(R^{1/\theta})}$ and base width $n_2^{-l_3(R^{1/\theta})}$ are somewhat larger than $R^{1/\theta}/n_1,R^{1/\theta}/n_2$ respectively since $l_3(R^{1/\theta})$ $\leq l_2(R^{1/\theta})\leq l_1(R^{1/\theta})$. We now iterate inside each cuboid constructed in \textbf{Step $3^{\circ}$} until we obtain a collection of cuboids each with base length $n_1^{-l_2(R^{1/\theta})}$ and base width $n_2^{-l_2(R^{1/\theta})}$, which takes further $l_2(R^{1/\theta})-l_3(R^{1/\theta})$ iterations. Note that we only need a factor of $|\pi_2\mathcal{D}|$ more covering sets since we may cover cuboids obtained in \textbf{Step $3^{\circ}$} simultaneously provided they share the same base.

       \textbf{Step $5^{\circ}$:} After \textbf{Step $4^{\circ}$}, we have a family of cuboids each with height and base width approximately $R^{1/\theta}$, whereas its base length $n_1^{-l_2(R^{1/\theta})}$ is still larger than $R^{1/\theta}/n_1$ since $l_2(R^{1/\theta})\leq l_1(R^{1/\theta})$. We continue to iterate inside each cuboid constructed in \textbf{Step $4^{\circ}$} until we obtain a collection of cuboids each with base length $n_1^{-l_1(R^{1/\theta})}$ (and thus approximately $R^{1/\theta}$), which takes further $l_1(R^{1/\theta})-l_2(R^{1/\theta})$ iterations. Note that we only need a factor of $|\pi_1\mathcal{D}|$ more covering sets since we may group together cuboids in \textbf{Step $4^{\circ}$} according to the first coordinate.

       Combining Steps $1^{\circ}-5^{\circ}$ and using \textbf{Notation 1}, we obtain that
       \begin{align*}
               &N(\Pi(Q(\omega,R)),R^{1/\theta})\\
               & \asymp \left(\prod\limits_{l=l_3(R)+1}^{l_2(R)} N(i_{l,1},i_{l,2})\right) \left(\prod\limits_{l=l_2(R)+1}^{l_1(R)} \sum_{v_{l,2}: (i_{l,1},v_{l,2})\in\pi_2 \mathcal{D}} N(i_{l,1},v_{l,2})\right) \\
               &\cdot\left(|\mathcal{D}|^{l_3(R^{1/\theta})-l_1(R)}\right) \left(|\pi_2\mathcal{D}|^{l_2(R^{1/\theta})-l_3(R^{1/\theta})}\right) \left(|\pi_1\mathcal{D}|^{l_1(R^{1/\theta})-l_2(R^{1/\theta})}\right)\\
               &\leq N_{\max}^{l_2(R)-l_3(R)} C_{\max}^{l_1(R)-l_2(R)} \left(|\mathcal{D}|^{l_3(R^{1/\theta})-l_1(R)}\right) \\
               &\cdot\left(|\pi_2\mathcal{D}|^{l_2(R^{1/\theta})-l_3(R^{1/\theta})}\right) \left(|\pi_1\mathcal{D}|^{l_1(R^{1/\theta})-l_2(R^{1/\theta})}\right) \\
               & \asymp N_{\max}^{\log R/\log n_3-\log R/\log n_2} \cdot C_{\max}^{\log R/\log n_2-\log R/\log n_1} \cdot |\mathcal{D}|^{\log R/\log n_1-\log R/\theta\log n_3} \\
               &\cdot |\pi_2\mathcal{D}|^{\log R/\theta\log n_3-\log R/\theta\log n_2} \cdot |\pi_1\mathcal{D}|^{\log R/\theta\log n_2-\log R/\theta\log n_1} \\
               & = R^{\log N_{\max}/\log n_3+\log (C_{\max}/N_{\max})/\log n_2+\log (|\mathcal{D}|/C_{\max})/\log n_1}\\
               &\cdot R^{\log (|\pi_2\mathcal{D}|/|\mathcal{D}|)/\theta\log n_3+\log (|\pi_1\mathcal{D}|/|\pi_2\mathcal{D}|)/\theta\log n_2-\log |\pi_1\mathcal{D}|/\theta\log n_1}.
       \end{align*}
       Taking logarithms and dividing by $(1-1/\theta)\log R$, we obtain that
       \begin{align*}
               &\frac{\log N(\Pi(Q(\omega,R)),R^{1/\theta})}{(1-1/\theta)\log R}\\
               &\leq \frac{\frac{\log N_{\max}}{\log n_3}+\frac{\log (C_{\max}/N_{\max})}{\log n_2}+\frac{\log (|\mathcal{D}|/C_{\max})}{\log n_1}}{1-1/\theta}\\
               &+\frac{\frac{\log (|\pi_2\mathcal{D}|/|\mathcal{D}|)}{\theta\log n_3}+\frac{\log (|\pi_1\mathcal{D}|/|\pi_2\mathcal{D}|)}{\theta\log n_2}-\frac{\log |\pi_1\mathcal{D}|}{\theta\log n_1}}{1-1/\theta}
               +\frac{O(1)}{\log R}\\
               &=\frac{\frac{\log|\pi_1\mathcal{D}|}{\log n_1}+\frac{\log (|\pi_2\mathcal{D}|/|\pi_1\mathcal{D}|)}{\log n_2}+\frac{\log (|\mathcal{D}|/|\pi_2\mathcal{D}|)}{\log n_3}}{1-\theta} \\
               &-\frac{\theta\left(\frac{\log (|\mathcal{D}|/C_{\max})}{\log n_1}+\frac{\log (C_{\max}/N_{\max})}{\log n_2}+\frac{\log N_{\max}}{\log n_3}\right)}{1-\theta}+\frac{O(1)}{\log R}.
       \end{align*}
       
       Letting $R\to 0$, this implies the desired upper bound. Observe that if we choose $\omega\in\mathcal{D}^{\infty}$ satisfying $\sum_{v_{l,2}: (i_{l,1},v_{l,2})\in\pi_2 \mathcal{D}} N(i_{l,1},v_{l,2})=C_{\max}$ for all $l_2(R)+1\leq l\leq l_1(R)$ and $N(i_{l,1},i_{l,2})=N_{\max}$ for all $l_3(R)+1\leq l\leq l_2(R)$, then the unique inequality in the above argument is replaced by equality and our covering estimates were optimal up to multiplicative constants. Therefore, in this way we will obtain the required lower bound for the Assouad spectrum.

       We now turn to calculate the formula of the lower spectrum. Repeating the above covering argument, we obtain that
       \begin{align*}
               &N(\Pi(Q(\omega,R)),R^{1/\theta})\\
               & \asymp \left(\prod\limits_{l=l_3(R)+1}^{l_2(R)} N(i_{l,1},i_{l,2})\right) \left(\prod\limits_{l=l_2(R)+1}^{l_1(R)} \sum_{v_{l,2}: (i_{l,1},v_{l,2})\in\pi_2 \mathcal{D}} N(i_{l,1},v_{l,2})\right)  \\
               &\cdot\left(|\mathcal{D}|^{l_3(R^{1/\theta})-l_1(R)}\right)\left(|\pi_2\mathcal{D}|^{l_2(R^{1/\theta})-l_3(R^{1/\theta})}\right) \left(|\pi_1\mathcal{D}|^{l_1(R^{1/\theta})-l_2(R^{1/\theta})}\right)\\
       \end{align*}
       whereas we proceed with uniform lower bounds instead:
       \begin{align*}
              &N(\Pi(Q(\omega,R)),R^{1/\theta})\\
              & \geq N_{\min}^{l_2(R)-l_3(R)} C_{\min}^{l_1(R)-l_2(R)} \left(|\mathcal{D}|^{l_3(R^{1/\theta})-l_1(R)}\right) \\
               &\cdot\left(|\pi_2\mathcal{D}|^{l_2(R^{1/\theta})-l_3(R^{1/\theta})}\right) \left(|\pi_1\mathcal{D}|^{l_1(R^{1/\theta})-l_2(R^{1/\theta})}\right) \\
               & \asymp N_{\min}^{\log R/\log n_3-\log R/\log n_2} \cdot C_{\min}^{\log R/\log n_2-\log R/\log n_1} \cdot |\mathcal{D}|^{\log R/\log n_1-\log R/\theta\log n_3} \\
               &\cdot |\pi_2\mathcal{D}|^{\log R/\theta\log n_3-\log R/\theta\log n_2} \cdot |\pi_1\mathcal{D}|^{\log R/\theta\log n_2-\log R/\theta\log n_1} \\
               & = R^{\log N_{\min}/\log n_3+\log (C_{\min}/N_{\min})/\log n_2+\log (|\mathcal{D}|/C_{\min})/\log n_1}\\
               &\cdot R^{\log (|\pi_2\mathcal{D}|/|\mathcal{D}|)/\theta\log n_3+\log (|\pi_1\mathcal{D}|/|\pi_2\mathcal{D}|)/\theta\log n_2-\log |\pi_1\mathcal{D}|/\theta\log n_1}.
       \end{align*}
       Taking logarithms and dividing by $(1-1/\theta)\log R$, we obtain that
       \begin{equation*}
           \begin{split}
               &\frac{\log N(\Pi(Q(\omega,R)),R^{1/\theta})}{(1-1/\theta)\log R} \\
               &\geq\frac{\frac{\log N_{\min}}{\log n_3}+\frac{\log (C_{\min}/N_{\min})}{\log n_2}+\frac{\log (|\mathcal{D}|/C_{\min})}{\log n_1}}{1-1/\theta} \\
               &+\frac{\frac{\log (|\pi_2\mathcal{D}|/|\mathcal{D}|)}{\theta\log n_3}+\frac{\log (|\pi_1\mathcal{D}|/|\pi_2\mathcal{D}|)}{\theta\log n_2}-\frac{\log |\pi_1\mathcal{D}|}{\theta\log n_1}}{1-1/\theta}+\frac{O(1)}{\log R}\\
               &=\frac{\frac{\log|\pi_1\mathcal{D}|}{\log n_1}+\frac{\log (|\pi_2\mathcal{D}|/|\pi_1\mathcal{D}|)}{\log n_2}+\frac{\log (|\mathcal{D}|/|\pi_2\mathcal{D}|)}{\log n_3}}{1-\theta}\\
               &-\frac{\theta\left(\frac{\log (|\mathcal{D}|/C_{\min})}{\log n_1}+\frac{\log (C_{\min}/N_{\min})}{\log n_2}+\frac{\log N_{\min}}{\log n_3}\right)}{1-\theta}+\frac{O(1)}{\log R}.
           \end{split}
       \end{equation*}
       Letting $R\to 0$, this implies the desired lower bound. Observe that if we choose $\omega\in\mathcal{D}^{\infty}$ satisfying $\sum_{v_{l,2}: (i_{l,1},v_{l,2})\in\pi_2 \mathcal{D}} N(i_{l,1},v_{l,2})=C_{\min}$ for all $l_2(R)+1\leq l\leq l_1(R)$ and $N(i_{l,1},i_{l,2})=N_{\min}$ for all $l_3(R)+1\leq l\leq l_2(R)$, then the unique inequality in the above argument is replaced by equality and our covering estimates were optimal up to multiplicative constants. Therefore, in this way we will obtain the required upper bound for the lower spectrum.
       \item if $\theta\in(\log n_1/\log n_3,\log n_1/\log n_2]$: we proceed with the covering argument in $(1)$, whereas this time we will obtain a collection of cuboids each with height approximately $R^{1/\theta}$ after $l_3(R^{1/\theta})$ steps, which is before the base length becomes smaller than $R$, since in this situation $l_3(R^{1/\theta})\leq l_1(R)$. For this reason, the third term (concerning powers of $|\mathcal{D}|$) is not required anymore. We have
       \begin{align*}
               &N(\Pi(Q(\omega,R)),R^{1/\theta})\\
               & \asymp \left(\prod\limits_{l=l_3(R)+1}^{l_2(R)} N(i_{l,1},i_{l,2})\right) \left(\prod\limits_{l=l_2(R)+1}^{l_3(R^{1/\theta})} \sum_{v_{l,2}: (i_{l,1},v_{l,2})\in\pi_2 \mathcal{D}} N(i_{l,1},v_{l,2})\right)   \\
               &\cdot\left(\prod\limits_{l=l_3(R^{1/\theta})+1}^{l_1(R)} N(i_{l,1})\right)\left(|\pi_2\mathcal{D}|^{l_2(R^{1/\theta})-l_1(R)}\right) \left(|\pi_1\mathcal{D}|^{l_1(R^{1/\theta})-l_2(R^{1/\theta})}\right)\\
               &\leq N_{\max}^{l_2(R)-l_3(R)} C_{\max}^{l_3(R^{1/\theta})-l_2(R)} M_{\max}^{l_1(R)-l_3(R^{1/\theta})} \\
               &\cdot\left(|\pi_2\mathcal{D}|^{l_2(R^{1/\theta})-l_1(R)}\right) \left(|\pi_1\mathcal{D}|^{l_1(R^{1/\theta})-l_2(R^{1/\theta})}\right) \\
               & \asymp N_{\max}^{\log R/\log n_3-\log R/\log n_2} \cdot C_{\max}^{\log R/\log n_2-\log R/\theta\log n_3} \cdot M_{\max}^{\log R/\theta\log n_3-\log R/\log n_1} \\
               &\cdot |\pi_2\mathcal{D}|^{\log R/\log n_1-\log R/\theta\log n_2} \cdot |\pi_1\mathcal{D}|^{\log R/\theta\log n_2-\log R/\theta\log n_1} \\
               & = R^{\log N_{\max}/\log n_3+\log (C_{\max}/N_{\max})/\log n_2+\log (|\pi_2\mathcal{D}|/M_{\max})/\log n_1}\\
               &\cdot R^{\log (M_{\max}/C_{\max})/\theta\log n_3+\log (|\pi_1\mathcal{D}|/|\pi_2\mathcal{D}|)/\theta\log n_2-\log |\pi_1\mathcal{D}|/\theta\log n_1}.
       \end{align*}
       Taking logarithms and dividing by $(1-1/\theta)\log R$, we obtain that
       \begin{equation*}
           \begin{split}
               &\frac{\log N(\Pi(Q(\omega,R)),R^{1/\theta})}{(1-1/\theta)\log R}  \\
               &\leq\frac{\frac{\log N_{\max}}{\log n_3}+\frac{\log (C_{\max}/N_{\max})}{\log n_2}+\frac{\log (|\pi_2\mathcal{D}|/M_{\max})}{\log n_1}}{1-1/\theta}\\
               &+\frac{\frac{\log (M_{\max}/C_{\max})}{\theta\log n_3}+\frac{\log (|\pi_1\mathcal{D}|/|\pi_2\mathcal{D}|)}{\theta\log n_2}-\frac{\log |\pi_1\mathcal{D}|}{\theta\log n_1}}{1-1/\theta}+\frac{O(1)}{\log R}\\
               &=\frac{\frac{\log|\pi_1\mathcal{D}|}{\log n_1}+\frac{\log (|\pi_2\mathcal{D}|/|\pi_1\mathcal{D}|)}{\log n_2}+\frac{\log(C_{\max}/M_{\max})}{\log n_3}}{1-\theta} \\
               &-\frac{\theta\left(\frac{\log (|\pi_2\mathcal{D}|/M_{\max})}{\log n_1}+\frac{\log (C_{\max}/N_{\max})}{\log n_2}+\frac{\log N_{\max}}{\log n_3}\right)}{1-\theta}+\frac{O(1)}{\log R}.
           \end{split}
       \end{equation*}
       Letting $R\to 0$, this implies the desired upper bound. Once again, one obtains the required lower bound by choosing $\omega\in\mathcal{D}^{\infty}$ satisfying $N(i_{l,1},i_{l,2})=N_{\max}$ for all $l_3(R)+1\leq l\leq l_2(R)$, $\sum_{v_{l,2}: (i_{l,1},v_{l,2})\in\pi_2 \mathcal{D}} N(i_{l,1},v_{l,2})=C_{\max}$ for all $l_2(R)+1\leq l\leq l_3(R^{1/\theta})$ and $N(i_{l,1})=M_{\max}$ for all $l_3(R^{1/\theta})+1\leq l\leq l_1(R)$. Similar to the argument in $(1)$, one can also obtain the formula of the lower spectrum.
       \item if $\theta\in(\log n_1/\log n_2,\log n_2/\log n_3]$, we proceed with the covering argument in $(1)$, whereas this time we will obtain a collection of cuboids each with both base width and height approximately $R^{1/\theta}$ after $l_2(R^{1/\theta})$ steps, which is before the base length becomes smaller than $R$, since in this situation $l_3(R^{1/\theta})\leq l_2(R^{1/\theta})\leq l_1(R)$. For this reason, the third term (concerning powers of $|\mathcal{D}|$) and the fourth term (concerning powers of $|\pi_2\mathcal{D}|$) are not required anymore. We have
       \begin{align*}
               &N(\Pi(Q(\omega,R)),R^{1/\theta})\\
               & \asymp \left(\prod\limits_{l=l_3(R)+1}^{l_2(R)} N(i_{l,1},i_{l,2})\right) \left(\prod\limits_{l=l_2(R)+1}^{l_3(R^{1/\theta})} \sum_{v_{l,2}: (i_{l,1},v_{l,2})\in\pi_2 \mathcal{D}} N(i_{l,1},v_{l,2})\right)\\
               &\cdot \left(\prod\limits_{l=l_3(R^{1/\theta})+1}^{l_2(R^{1/\theta})} N(i_{l,1})\right) \cdot \left(|\pi_1\mathcal{D}|^{l_1(R^{1/\theta})-l_1(R)}\right)\\
               &\leq N_{\max}^{l_2(R)-l_3(R)} C_{\max}^{l_3(R^{1/\theta})-l_2(R)} M_{\max}^{l_2(R^{1/\theta})-l_3(R^{1/\theta})} 
               \cdot \left(|\pi_1\mathcal{D}|^{l_1(R^{1/\theta})-l_1(R)}\right) \\
               & \asymp N_{\max}^{\log R/\log n_3-\log R/\log n_2} \cdot C_{\max}^{\log R/\log n_2-\log R/\theta\log n_3} \\
               &\cdot M_{\max}^{\log R/\theta\log n_3-\log R/\theta\log n_2} \cdot |\pi_1\mathcal{D}|^{\log R/\log n_1-\log R/\theta\log n_1} \\
               & = R^{\log N_{\max}/\log n_3+\log (C_{\max}/N_{\max})/\log n_2+\log |\pi_1\mathcal{D}|/\log n_1}\\
               &\cdot R^{\log (M_{\max}/C_{\max})/\theta\log n_3-\log M_{\max}/\theta\log n_2-\log |\pi_1\mathcal{D}|/\theta\log n_1}.
       \end{align*}
       Taking logarithms and dividing by $(1-1/\theta)\log R$, we obtain that
       \begin{equation*}
           \begin{split}
               &\frac{\log N(\Pi(Q(\omega,R)),R^{1/\theta})}{(1-1/\theta)\log R}\\
               &\leq \frac{\frac{\log N_{\max}}{\log n_3}+\frac{\log (C_{\max}/N_{\max})}{\log n_2}+\frac{\log |\pi_1\mathcal{D}|}{\log n_1}}{1-1/\theta}\\
               &+\frac{\frac{\log (M_{\max}/C_{\max})}{\theta\log n_3}-\frac{\log M_{\max}}{\theta\log n_2}-\frac{\log |\pi_1\mathcal{D}|}{\theta\log n_1}}{1-1/\theta}+\frac{O(1)}{\log R}\\
               &=\frac{\log|\pi_1\mathcal{D}|}{\log n_1}+\frac{\frac{\log M_{\max}}{\log n_2}+\frac{\log (C_{\max}/M_{\max})}{\log n_3}-\theta\left(\frac{\log (C_{\max}/N_{\max})}{\log n_2}+\frac{\log N_{\max}}{\log n_3}\right)}{1-\theta} +\frac{O(1)}{\log R}.
           \end{split}
       \end{equation*}
       Letting $R\to 0$, this implies the desired upper bound. Once more, one obtains the required lower bound by choosing $\omega\in\mathcal{D}^{\infty}$ satisfying $N(i_{l,1},i_{l,2})=N_{\max}$ for all $l_3(R)+1\leq l\leq l_2(R)$, $\sum_{v_{l,2}: (i_{l,1},v_{l,2})\in\pi_2 \mathcal{D}} N(i_{l,1},v_{l,2})=C_{\max}$ for all $l_2(R)+1\leq l\leq l_3(R^{1/\theta})$ and $N(i_{l,1})=M_{\max}$ for all $l_3(R^{1/\theta})+1\leq l\leq l_2(R^{1/\theta})$. Similar to the argument in $(1)$, one can also obtain the formula of the lower spectrum.

    \item if $\theta\in(\log n_2/\log n_3,1)$, we proceed with the covering argument in $(1)$, whereas this time we will obtain a collection of middle sized cuboids each with height approximately $R^{1/\theta}$ after $l_3(R^{1/\theta})$ steps, which is before the base width becomes smaller than $R$, since $l_3(R^{1/\theta})\leq l_2(R)$. Moreover, after further $l_2(R^{1/\theta})-l_3(R^{1/\theta})$ steps, we will obtain a large family of smaller sized cuboids each with base width approximately $R^{1/\theta}$, which is before the base length becomes smaller than $R$, since $l_2(R^{1/\theta})\leq l_1(R)$. For this reason, the third term (concerning powers of $|\mathcal{D}|$) and the fourth term (concerning powers of $|\pi_2\mathcal{D}|$) disappear. We have
       \begin{equation*}
           \begin{split}
               &N(\Pi(Q(\omega,R)),R^{1/\theta})\\
               & \asymp \left(\prod\limits_{l=l_3(R)+1}^{l_3(R^{1/\theta})} N(i_{l,1},i_{l,2})\right) \left(\prod\limits_{l=l_2(R)+1}^{l_2(R^{1/\theta})} N(i_{l,1})\right) \cdot \left(|\pi_1\mathcal{D}|^{l_1(R^{1/\theta})-l_1(R)}\right)\\
               &\leq N_{\max}^{l_3(R^{1/\theta})-l_3(R)} M_{\max}^{l_2(R^{1/\theta})-l_2(R)} 
               \cdot \left(|\pi_1\mathcal{D}|^{l_1(R^{1/\theta})-l_1(R)}\right) \\
               & \asymp N_{\max}^{\log R/\log n_3-\log R/\theta\log n_3} \cdot M_{\max}^{\log R/\log n_2-\log R/\theta\log n_2}\\
               &\cdot |\pi_1\mathcal{D}|^{\log R/\log n_1-\log R/\theta\log n_1} \\
               & = R^{(1-1/\theta)(\log N_{\max}/\log n_3+\log M_{\max}/\log n_2+\log |\pi_1\mathcal{D}|/\log n_1}).
           \end{split}
       \end{equation*}
       Taking logarithms and dividing by $(1-1/\theta)\log R$, we obtain that
       \begin{equation*}
               \frac{\log N(\Pi(Q(\omega,R)),R^{1/\theta})}{(1-1/\theta)\log R} \leq \frac{\log |\pi_1\mathcal{D}|}{\log n_1}+\frac{\log M_{\max}}{\log n_2}+\frac{\log N_{\max}}{\log n_3} +\frac{O(1)}{\log R}.
       \end{equation*}
       Letting $R\to 0$, this implies the desired upper bound. Once more, one obtains the required lower bound by choosing $\omega\in\mathcal{D}^{\infty}$ satisfying $N(i_{l,1},i_{l,2})=N_{\max}$ for all $l_3(R)+1\leq l\leq l_3(R^{1/\theta})$ and $N(i_{l,1})=M_{\max}$ for all $l_2(R)+1\leq l\leq l_2(R^{1/\theta})$. Similar to the argument in $(1)$, one can also obtain the formula of the lower spectrum.  Q.E.D.
    \end{enumerate}

\subsection{Proof of Theorem \ref{spectra case(2)}}
    Recall that we are assuming $\log n_2/\log n_3<\log n_1/\log n_2$. Observe that in cases $0<\theta\leq\log n_1/\log n_3, \log n_1/\log n_3<\theta\leq \log n_2/\log n_3, \log n_1/\log n_2<\theta<1$, the Assouad and lower spectra have the same forms as in cases $0<\theta\leq\log n_1/\log n_3$, $\log n_1/\log n_3<\theta\leq \log n_1/\log n_2, \log n_2/\log n_3<\theta<1$ of Theorem \ref{spectra formula} respectively, which is ensured by Lemma \ref{scale comparison}. Thus it suffices to consider the case $\log n_2/\log n_3<\theta\leq \log n_1/\log n_2$. We proceed with the covering argument in Theorem \ref{spectra formula} $(1)$, whereas this time we will obtain a collection of cuboids each with height approximately $R^{1/\theta}$ after $l_3(R^{1/\theta})$ steps, which is before the base width becomes smaller than $R$, since in this situation $l_3(R^{1/\theta})\leq l_2(R)$. For this reason, the third term (concerning powers of $|\mathcal{D}|$) is not required anymore. 
    
    We have
       \begin{align*}
               &N(\Pi(Q(\omega,R)),R^{1/\theta})\\
               & \asymp \left(\prod\limits_{l=l_3(R)+1}^{l_3(R^{1/\theta})} N(i_{l,1},i_{l,2})\right) \left(\prod\limits_{l=l_2(R)+1}^{l_1(R)} N(i_{l,1})\right)  \\
               &\cdot\left(|\pi_2\mathcal{D}|^{l_2(R^{1/\theta})-l_1(R)}\right) \left(|\pi_1\mathcal{D}|^{l_1(R^{1/\theta})-l_2(R^{1/\theta})}\right)\\
               &\leq N_{\max}^{l_3(R^{1/\theta})-l_3(R)} M_{\max}^{l_1(R)-l_2(R)} 
               \cdot\left(|\pi_2\mathcal{D}|^{l_2(R^{1/\theta})-l_1(R)}\right) \left(|\pi_1\mathcal{D}|^{l_1(R^{1/\theta})-l_2(R^{1/\theta})}\right) \\
               & \asymp N_{\max}^{\log R/\log n_3-\log R/\theta\log n_3} \cdot M_{\max}^{\log R/\log n_2-\log R/\log n_1} \\
               &\cdot |\pi_2\mathcal{D}|^{\log R/\log n_1-\log R/\theta\log n_2} \cdot |\pi_1\mathcal{D}|^{\log R/\theta\log n_2-\log R/\theta\log n_1} \\
               & = R^{\log N_{\max}/\log n_3+\log M_{\max}/\log n_2+\log (|\pi_2\mathcal{D}|/M_{\max})/\log n_1}\\
               &\cdot R^{-\log N_{\max}/\theta\log n_3+\log (|\pi_1\mathcal{D}|/|\pi_2\mathcal{D}|)/\theta\log n_2-\log |\pi_1\mathcal{D}|/\theta\log n_1}.
       \end{align*}
       Taking logarithms and dividing by $(1-1/\theta)\log R$, we obtain that
       \begin{equation*}
           \begin{split}
              & \frac{\log N(\Pi(Q(\omega,R)),R^{1/\theta})}{(1-1/\theta)\log R} \\
              &\leq \frac{\frac{\log N_{\max}}{\log n_3}+\frac{\log M_{\max}}{\log n_2}+\frac{\log (|\pi_2\mathcal{D}|/M_{\max})}{\log n_1}}{1-1/\theta} \\
              &+\frac{-\frac{\log N_{\max}}{\theta\log n_3}+\frac{\log (|\pi_1\mathcal{D}|/|\pi_2\mathcal{D}|)}{\theta\log n_2}-\frac{\log |\pi_1\mathcal{D}|}{\theta\log n_1}}{1-1/\theta}+\frac{O(1)}{\log R}\\
               &=\frac{\frac{\log|\pi_1\mathcal{D}|}{\log n_1}+\frac{\log (|\pi_2\mathcal{D}|/|\pi_1\mathcal{D}|)}{\log n_2}-\theta\left(\frac{\log (|\pi_2\mathcal{D}|/M_{\max})}{\log n_1}+\frac{\log M_{\max}}{\log n_2}\right)}{1-\theta}+\frac{\log N_{\max}}{\log n_3} +\frac{O(1)}{\log R}.
           \end{split}
       \end{equation*}
       Letting $R\to 0$, this implies the desired upper bound. Once again, one obtains the required lower bound by choosing $\omega\in\mathcal{D}^{\infty}$ satisfying $N(i_{l,1},i_{l,2})=N_{\max}$ for all $l_3(R)+1\leq l\leq l_3(R^{1/\theta})$ and $N(i_{l,1})=M_{\max}$ for all $l_2(R)+1\leq l\leq l_1(R)$. Similar to the argument in Theorem \ref{spectra formula}, one can also obtain the formula of the lower spectrum.  Q.E.D.

\subsection{Proof of Theorem \ref{spectra case(3)}}
    We may obtain the formulae for the Assouad and lower spectra immediately from Theorem \ref{spectra formula}, due to Lemma \ref{scale comparison}.

\subsection{Proof of Proposition \ref{Lipschitz equivalence}}
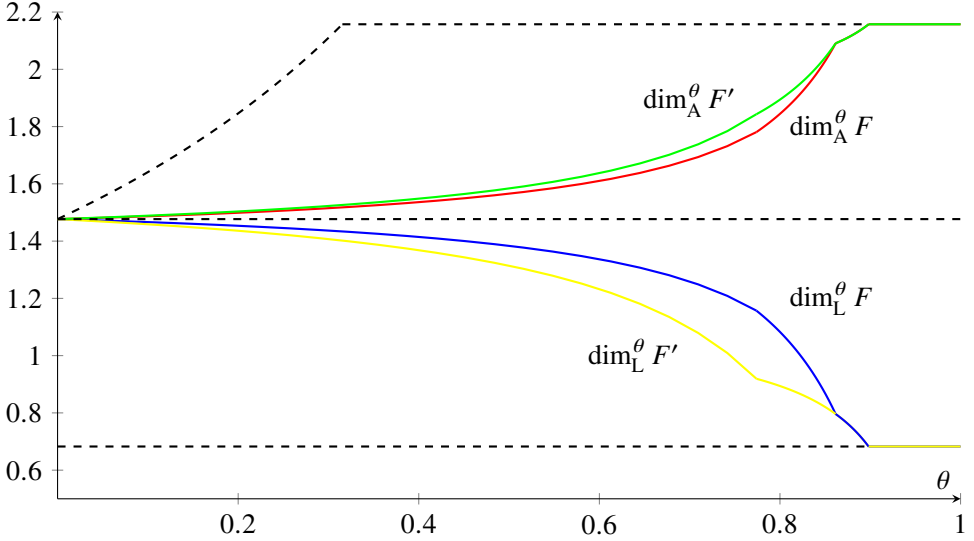
\begin{figure}[t]
    \caption{Plots of the Assouad and lower spectra for the two sponges constructed in the proof of Proposition \ref{Lipschitz equivalence}. The phase transition for $\adim^{\theta} F, \adim^{\theta} F'$ and $\ladim^{\theta} F$ occurs at $\theta=\log5/\log8,\log6/\log8$ and $\log5/\log6$, and the phase transition for $\ladim^{\theta} F'$ occurs at $\theta=\log5/\log8$ and $\log5/\log6$. The general bounds from \cite[Proposition 3.1]{FY18} and \cite[Proposition 3.9]{FY18} are drawn with dashed lines. Also, the unique phase transition in the general bound occurs at $\theta=1-\bdim F/\adim F\approx 0.315$.}
    \label{Plots of two sponges}
      \begin{tikzpicture}
        \begin{axis}[
    axis lines = middle,
    xlabel = $\theta$,
    ymin=0.5, ymax=2.2,
    xmin=0, xmax=1,
    width=13.5cm,
    height=8cm,
    xtick distance=0.2,
    ]
 \addplot [
     domain=0:log10(5)/log10(8), 
     color=blue,
     thick,
 ]
 {(log10(3)/log10(5) + log10(2)/log10(6) + log10(7/3)/log10(8) - x*((log10(14/3))/(log10(5)) +  log10(3)/log10(6)))/(1-x)};

 \addplot [
     domain=0:log10(5)/log10(8), 
     color=yellow,
     thick,
 ]
 {(log10(3)/log10(5) + log10(2)/log10(6) + log10(7/3)/log10(8) - x*((log10(14))/(log10(5)) ))/(1-x)};

  \addplot [
     domain=0:log10(5)/log10(8), 
     color=red,
     thick,
 ]
 {(log10(3)/log10(5) + log10(2)/log10(6) + log10(7/3)/log10(8) - x*((log10(7/3))/(log10(5)) + log10(6)/log10(8)))/(1-x)};

 \addplot [
     domain=0:log10(5)/log10(8), 
     color=green,
     thick,
 ]
 {(log10(3)/log10(5) + log10(2)/log10(6) + log10(7/3)/log10(8) - x*((log10(7/4))/(log10(5)) + (log10(4/3))/(log10(6)) + log10(6)/log10(8)))/(1-x)};


  \addplot [
     domain=log10(5)/log10(8):log10(6)/log10(8),
     color=blue,
     thick,
 ]
 {(log10(3)/log10(5) + log10(2)/log10(6) + log10(3)/log10(8)  - x*(log10(6)/log10(5) + log10(3)/log10(6)))/(1-x) };

 \addplot [
     domain=log10(5)/log10(8):log10(5)/log10(6),
     color=yellow,
     thick,
 ]
 {(log10(3)/log10(5) + log10(2)/log10(6) - x*((log10(6))/(log10(5)) ))/(1-x)};

  \addplot [
     domain=log10(5)/log10(8):log10(6)/log10(8),
     color=red,
     thick,
 ]
 {(log10(3)/log10(5) + log10(2)/log10(6) + log10(2)/log10(8)  - x*(log10(2)/log10(5) + log10(6)/log10(8)))/(1-x) };

 \addplot [
     domain=log10(5)/log10(8):log10(6)/log10(8), 
     color=green,
     thick,
 ]
 {(log10(3)/log10(5) + log10(2)/log10(6) + log10(8/3)/log10(8) - x*((log10(2))/(log10(5)) + (log10(4/3))/(log10(6)) + log10(6)/log10(8)))/(1-x)};

\addplot [
    domain=log10(6)/log10(8):log10(5)/log10(6),
    color=blue,
    thick,
]
{(log10(3)/log10(5) + log10(2)/log10(6) - x*(log10(6)/log10(5)))/(1-x) };

\addplot [
    domain=log10(6)/log10(8):log10(5)/log10(6),
    color=red,
    thick,
]
{(log10(3)/log10(5) + log10(2)/log10(6) + log10(6)/log10(8)  - x*(log10(2)/log10(5) + log10(3)/log10(6) + log10(6)/log10(8)))/(1-x) };

 \addplot [
      domain=log10(6)/log10(8):log10(5)/log10(6), 
      color=green,
      thick,
  ]
  {(log10(3)/log10(5) + log10(2)/log10(6) + log10(6)/log10(8) - x*((log10(2))/(log10(5)) + (log10(3))/(log10(6)) + log10(6)/log10(8)))/(1-x)};

\addplot [
    domain=log10(5)/log10(6):1, 
    color=blue,
    thick,
]
{log10(3)/log10(5)};

\addplot [
    domain=log10(5)/log10(6):1, 
    color=yellow,
    thick,
]
{log10(3)/log10(5)};

\addplot [
    domain=log10(5)/log10(6):1, 
    color=red,
    thick,
]
{log10(3)/log10(5) + log10(3)/log10(6) + log10(6)/log10(8) };

\addplot [
    domain=log10(5)/log10(6):1, 
    color=green,
    thick,
]
{log10(3)/log10(5) + log10(3)/log10(6) + log10(6)/log10(8)};

\addplot [
    domain=0:log10(5)/log10(6),
    dashed,
    thick,
]
{log10(3)/log10(5)};

\addplot [
    domain=0:(log10(3/2)/log10(6)+log10(18/7)/log10(8))/(log10(3)/log10(5)+log10(3)/log10(6)+log10(6)/log10(8)), 
    dashed,
    thick,
]
{(log10(3)/log10(5) + log10(2)/log10(6) + log10(7/3)/log10(8) )/(1-x)};

\addplot [
    domain=(log10(3/2)/log10(6)+log10(18/7)/log10(8))/(log10(3)/log10(5)+log10(3)/log10(6)+log10(6)/log10(8)):log10(5)/log10(6), 
    dashed,
    thick,
]
{log10(3)/log10(5) + log10(3)/log10(6) + log10(6)/log10(8)};

\addplot [
    domain=0:1, 
    dashed,
    thick,
]
{log10(3)/log10(5) + log10(2)/log10(6) + log10(7/3)/log10(8) };

\node [right] at (axis cs:0.8,1.2){$\ladim^{\theta} F$};
\node [left] at (axis cs:0.7,1){$\ladim^{\theta} F'$};
\node [above] at (axis cs:0.7,1.8){$\adim^{\theta} F'$};
\node [right] at (axis cs:0.8,1.8){$\adim^{\theta} F$};

        \end{axis}
      \end{tikzpicture}
      \hfill

\end{figure}

  Let $n_1=5,n_2=6$ and $n_3=8$.
  We construct the first sponge $F$ generated by
  \begin{equation*}
      \begin{split}
          \mathcal{D}=\{&(0,0,0),(0,0,1),(0,0,2),(0,0,3),(0,0,4),(0,0,5),(1,1,0),(1,1,2),\\
          &(1,1,3),(1,1,5),(1,3,4),(2,1,0),(2,2,2),(2,4,6)\}
      \end{split}
  \end{equation*}
  and the second sponge $F'$ generated by 
  \begin{equation*}
      \begin{split}
          \mathcal{D}'=\{&(0,0,0),(0,0,1),(0,0,2),(0,0,3),(0,0,4),(0,0,5),(0,2,2),(0,4,6),\\
          &(1,1,0),(1,1,2),(1,1,3),(1,1,5),(1,3,4),(2,1,0)\}.
      \end{split}
  \end{equation*}
  Then for both sponges one has 
  $|\mathcal{D}|=|\mathcal{D}'|=14,|\pi_2\mathcal{D}|=|\pi_2\mathcal{D}'|=6,|\pi_1\mathcal{D}|=|\pi_1\mathcal{D}'|=3,M_{\max}=M'_{\max}=3,M_{\min}=M'_{\min}=1,N_{\max}=N'_{\max}=6,N_{\min}=N'_{\min}=1$, and therefore by Theorem \ref{Kenyon-Peres-Fraser-Howroyd} one further has 
    \begin{equation*}
    \begin{split}
        &\adim F=\adim F'=\frac{\log 3}{\log 5}+\frac{\log 3}{\log 6}+\frac{\log 6}{\log 8}\approx 2.156,\\
        &\bdim F=\bdim F'=\frac{\log 3}{\log 5}+\frac{\log 2}{\log 6}+\frac{\log 7/3}{\log 8}\approx 1.477,\\
        &\ladim F=\ladim F'=\frac{\log 3}{\log 5}\approx 0.683.
    \end{split}
   \end{equation*}
 They are defined on the same grid and share the same box, lower and Assouad dimensions. However, they have different spectra since for the first sponge
 \begin{equation*}
     C_{\max}=N(0,0)=6, C_{\min}=N(2,1)+N(2,2)+N(2,4)=3
 \end{equation*}
 whereas for the second sponge 
 \begin{equation*}
   C'_{\max}=N(0,0)+N(0,2)+N(0,4)=8, C'_{\min}=N(2,1)=1.
 \end{equation*}
 Thus, the Bedford--McMullen sponges constructed above are not bi-Lipschitz equivalent since bi-Lipschitz maps preserve Assouad and lower spectra, see \cite[Corollary 4.9]{FY18}.
 Notice that $\log n_2/\log n_3<\log n_1/\log n_2$ for both sponges and we may apply Theorem \ref{spectra case(2)} to compute their spectra. For illustration see Figure \ref{Plots of two sponges}, where $\theta\mapsto\adim^{\theta} F$ is plotted in red, $\theta\mapsto\adim^{\theta} F'$ is plotted in green, $\theta\mapsto\ladim^{\theta} F$ is plotted in blue and $\theta\mapsto\ladim^{\theta} F'$ is plotted in yellow.  Q.E.D.

\subsection{Proof of Theorem \ref{general dimension}}
    Fix $\omega=(\bar\imath_1,\bar\imath_2,\ldots)=((i_{1,1},\ldots,i_{1,d}),(i_{2,1},\ldots,i_{2,d}),\ldots)\in\mathcal{D}^{\infty}$ and $R\in(0,1)$. Recall that $l_q(R)$ is the unique integer satisfying \eqref{lq} for $q=1,\ldots,d$ and $\Pi(Q(\omega,R))$ is a hypercuboid in $\R^d$ with all side lengths comparable to $R$.
    \begin{enumerate}[$(1)$]
        \item Note that $\log n_1/\log n_d$ is the minimal ratio among $\{\log n_p/\log n_q\}_{1\leq p<q\leq d}$ since we assume that the sequence $\{n_q\}_{1\leq q\leq d}$ is strictly increasing. If $\theta\in(0,\log n_1/\log n_d]$, direct calculation shows that 
        \begin{equation*}
            l_d(R)\leq l_{d-1}(R)\leq\ldots\leq l_1(R)\leq l_d(R^{1/\theta})\leq l_{d-1}(R^{1/\theta})\leq\ldots\leq l_1(R^{1/\theta}).
        \end{equation*}
         The approximate cube $\Pi(Q(\omega,R))$ is composed of 
         \begin{equation*}
             \begin{split}
                &\prod\limits_{l=l_d(R)+1}^{l_{d-1}(R)}N(i_{l,1},\ldots,i_{l,d-1})\\
                &\cdot\left(\prod\limits_{k=1}^{d-2}\prod\limits_{l=l_{k+1}(R)+1}^{l_k(R)}   \sum_{v_{l,k+1}:(i_{l,1},\ldots,i_{l,k},v_{l,k+1})\in\pi_{k+1}\mathcal{D}} \sum_{v_{l,k+2}:(i_{l,1},\ldots,i_{l,k},v_{l,k+1},v_{l,k+2})\in\pi_{k+2}\mathcal{D}}\right.\\
               &\cdots\left.\sum_{v_{l,d-1}:(i_{l,1},\ldots,i_{l,k},v_{l,k+1},\ldots,v_{l,d-1})\in\pi_{d-1}\mathcal{D}}N(i_{l,1},\ldots,i_{l,k},v_{l,k+1},\ldots,v_{l,d-1})\right)\\
             \end{split}
         \end{equation*}
         hypercuboids with all side lengths $n_q^{-l_1(R)}, q=1,\ldots,d$. We iterate inside each hypercuboid, i.e. we pick up $|\mathcal{D}|$ smaller hypercuboids, until the minimal side length is roughly $R^{1/\theta}$. This will take $l_d(R^{1/\theta})-l_1(R)$ iterations. We finally cover the resulting collection by small sets of diameter $R^{1/\theta}$.
         Putting these estimates together and using \textbf{Notation 2}, as in the proof of Theorem \ref{spectra formula}, we may obtain
       \begin{align*}
               &N(\Pi(Q(\omega,R)),R^{1/\theta})\\
               & \asymp \prod\limits_{l=l_d(R)+1}^{l_{d-1}(R)}N(i_{l,1},\ldots,i_{l,d-1})\\
               &\cdot\left(\prod\limits_{k=1}^{d-2}\prod\limits_{l=l_{k+1}(R)+1}^{l_k(R)}   \sum_{v_{l,k+1}:(i_{l,1},\ldots,i_{l,k},v_{l,k+1})\in\pi_{k+1}\mathcal{D}} \sum_{v_{l,k+2}:(i_{l,1},\ldots,i_{l,k},v_{l,k+1},v_{l,k+2})\in\pi_{k+2}\mathcal{D}}\right.\\
               &\cdots\left.\sum_{v_{l,d-1}:(i_{l,1},\ldots,i_{l,k},v_{l,k+1},\ldots,v_{l,d-1})\in\pi_{d-1}\mathcal{D}}N(i_{l,1},\ldots,i_{l,k},v_{l,k+1},\ldots,v_{l,d-1})\right)\\
               &\cdot\left(|\mathcal{D}|^{l_d(R^{1/\theta})-l_1(R)}\right) 
               \cdot\left(\prod\limits_{k=1}^{d-1}|\pi_k\mathcal{D}|^{l_k(R^{1/\theta})-l_{k+1}(R^{1/\theta})}\right)\\
               &\leq \prod\limits_{k=1}^{d-1}N_{k,\max}^{l_k(R)-l_{k+1}(R)} \left(|\mathcal{D}|^{l_d(R^{1/\theta})-l_1(R)}\right) \left(\prod\limits_{k=1}^{d-1}|\pi_k\mathcal{D}|^{l_k(R^{1/\theta})-l_{k+1}(R^{1/\theta})}\right)  \\
               & \asymp \prod\limits_{k=1}^{d-1} N_{k,\max}^{\log R/\log n_{k+1}-\log R/\log n_k} \cdot |\mathcal{D}|^{\log R/\log n_1-\log R/\theta\log n_d} \\
               &\cdot \prod\limits_{k=1}^{d-1}|\pi_k\mathcal{D}|^{\log R/\theta\log n_{k+1}-\log R/\theta\log n_k}  \\
               & = R^{\sum\limits_{k=1}^{d-1}(\log N_{k,\max}/\log n_{k+1}-\log N_{k,\max}/\log n_{k})+\log|\mathcal{D}|/\log n_1-\log|\mathcal{D}|/\theta\log n_d}\\
               &\cdot R^{\sum\limits_{k=1}^{d-1}(\log|\pi_k\mathcal{D}|/\theta\log n_{k+1}-\log|\pi_k\mathcal{D}|/\theta\log n_k)}.
       \end{align*}
       Taking logarithms and dividing by $(1-1/\theta)\log R$, we further obtain that
       \begin{align*}
               &\frac{\log N(\Pi(Q(\omega,R)),R^{1/\theta})}{(1-1/\theta)\log R} \\
               &\leq\frac{\frac{\log|\pi_1\mathcal{D}|}{\log n_1}+\sum\limits_{k=1}^{d-1}\frac{\log (|\pi_{k+1}\mathcal{D}|/|\pi_{k}\mathcal{D}|)}{\log n_{k+1}}}{1-\theta} \\
               &-\frac{\theta\left(\frac{\log (|\mathcal{D}|/N_{1,\max})}{\log n_1}+\sum\limits_{k=1}^{d-2}\frac{\log (N_{k,\max}/N_{k+1,\max})}{\log n_{k+1}}+\frac{\log N_{d-1,\max}}{\log n_d}\right)}{1-\theta}+\frac{O(1)}{\log R}.
       \end{align*}
       Letting $R\to 0$, this implies the desired upper bound. Observe that if we choose $\omega\in\mathcal{D}^{\infty}$ satisfying $N(i_{l,1},\ldots,i_{l,d-1})=N_{d-1,\max}$ for all $l_{d}(R)+1\leq l\leq l_{d-1}(R)$ and 
       \begin{equation*}
           \begin{split}
               &\sum_{v_{l,k+1}:(i_{l,1},\ldots,i_{l,k},v_{l,k+1})\in\pi_{k+1}\mathcal{D}} \sum_{v_{l,k+2}:(i_{l,1},\ldots,i_{l,k},v_{l,k+1},v_{l,k+2})\in\pi_{k+2}\mathcal{D}} \cdots\\
               &\sum_{v_{l,d-1}:(i_{l,1},\ldots,i_{l,k},v_{l,k+1},\ldots,v_{l,d-1})\in\pi_{d-1}\mathcal{D}}N(i_{l,1},\ldots,i_{l,k},v_{l,k+1},\ldots,v_{l,d-1})=N_{k,\max}
           \end{split}
       \end{equation*}
       for all $l_{k+1}(R)+1\leq l\leq l_k(R)$ and $1\leq k\leq d-2$,
       then the unique inequality in the above argument is replaced by equality and our covering estimates were optimal up to multiplicative constants. Therefore, in this way we will obtain the required lower bound for the Assouad spectrum. Similar to the argument in Theorem \ref{spectra formula}, one can also obtain the formula of the lower spectrum.
       \item If $\theta\in[\max_{1\leq p<q\leq d}\log n_p/\log n_q,1)$, then
       \begin{equation*}
            l_d(R)\leq l_d(R^{1/\theta})\leq l_{d-1}(R)\leq l_{d-1}(R^{1/\theta})\leq\ldots\leq l_1(R)\leq l_1(R^{1/\theta}).
        \end{equation*}
        Following the covering strategy in Theorem \ref{spectra formula}, we may obtain
         \begin{align*}
               &N(\Pi(Q(\omega,R)),R^{1/\theta})\\
               & \asymp \left(\prod\limits_{k=1}^{d-1}\prod\limits_{l=l_{k+1}(R)+1}^{l_{k+1}(R^{1/\theta})} N(i_{l,1},\ldots,i_{l,k})\right) \cdot \left(|\pi_1\mathcal{D}|^{l_1(R^{1/\theta})-l_1(R)}\right)\\
               &\leq N_{d-1,\max}^{l_{d}(R^{1/\theta})-l_{d}(R)}\prod\limits_{k=1}^{d-2}M_{k,\max}^{l_{k+1}(R^{1/\theta})-l_{k+1}(R)} \cdot \left(|\pi_1\mathcal{D}|^{l_1(R^{1/\theta})-l_1(R)}\right) \\
               & \asymp N_{d-1,\max}^{\log R/\log n_{d}-\log R/\theta\log n_{d}}\prod\limits_{k=1}^{d-2}M_{k,\max}^{\log R/\log n_{k+1}-\log R/\theta\log n_{k+1}} \\
               &\cdot |\pi_1\mathcal{D}|^{\log R/\log n_1-\log R/\theta\log n_1} \\
               & = R^{(1-1/\theta)(\log |\pi_1\mathcal{D}|/\log n_1+\sum\limits_{k=1}^{d-2}\log M_{k,\max}/\log n_{k+1}+N_{d-1,\max}/\log n_{d})}.
       \end{align*}
       Taking logarithms and dividing by $(1-1/\theta)\log R$, we obtain that
       \begin{equation*}
       \begin{split}
           &\frac{\log N(\Pi(Q(\omega,R)),R^{1/\theta})}{(1-1/\theta)\log R} \\
           &\leq \frac{\log |\pi_1\mathcal{D}|}{\log n_1}+\sum\limits_{k=1}^{d-2}\frac{\log M_{k,\max}}{\log n_{k+1}}+\frac{\log N_{d-1,\max}}{\log n_{d}} +\frac{O(1)}{\log R}.
       \end{split}
       \end{equation*}
       Letting $R\to 0$, this implies the desired upper bound. Once more, one obtains the required lower bound by choosing $\omega\in\mathcal{D}^{\infty}$ satisfying $N(i_{l,1},\ldots,i_{l,d-1})=N_{d-1,\max}$ and $N(i_{l,1},\ldots,i_{l,k})=M_{k,\max}$ for all $l_{k+1}(R)+1\leq l\leq l_{k+1}(R^{1/\theta})$ and $1\leq k\leq d-2$. Similar to the argument in Theorem \ref{spectra formula}, one can also obtain the formula of the lower spectrum.  Q.E.D.
    \end{enumerate}


\begin{ack}
The original way to define $C_{\max}$ and $C_{\min}$ in the proof of Theorem \ref{spectra formula} was incorrect. I would like to thank an anonymous reviewer for pointing this out and kindly suggesting the correct way (in the setup of $d=3$) to me, also for several useful comments which lead to improvements in the statement of the main results.
\end{ack}

\begin{funding}
The author was partially supported by the China Postdoctoral Science Foundation (2025M773065), the National Key R\&D Program of China (No. 2024YFA1013602, 2024YFA1013600) and NSFC grant (12090012, 12090010).
\end{funding}


\end{document}